\documentclass[12pt]{article}
\usepackage[latin1]{inputenc}
\usepackage{amssymb,amsmath,amsfonts}
\newtheorem{theorem}{Theorem}

\newtheorem{remark}{Remark}
\newtheorem{example}{Example}
\def\<{\langle}
\def\>{\rangle}
\newcommand{\R}{\mathbb{R}}
\newcommand{\Ric}{\mathrm{Ric}}
\newcommand{\po}{{\hspace*{-1ex}}{\bf .  }}
\newcommand\be{\begin{equation}} 
\newcommand\ee{\end{equation}}
\def\bea{\begin{eqnarray*} }
\def\eea{\end{eqnarray*} }
\def\proof{\noindent\emph{ Proof: }}
\begin{document}

\title{Entire unbounded constant mean curvature\\ Killing graphs}
\author{M. Dajczer and J. H. de Lira\thanks
{Partially supported by CNPq and FUNCAP.}}
\date{}
\maketitle

\begin{abstract} 
In this paper, we provided conditions for an entire constant mean 
curvature Killing graph lying inside a possible unbounded region
to be necessarily a slice.  
\end{abstract}

Let  $N^{n+1}$ denote a complete Riemannian manifold carrying a warped product 
structure $N^{n+1}=M^n\times_\rho\R$ with warping function $\rho\in C^\infty(M)$ 
such that  $M^n$ is noncompact. Let $s\in\R$ parametrize the line factor in 
$N^{n+1}$. Then $Y=\partial/\partial s$ is a Killing vector field free of 
singularities with integrable orthogonal distribution and $\rho=|Y|$.  
Moreover, the leaves $M^n\times\{s_0\},\, s_0\in\R$, of the orthogonal  
distribution to $Y$ form a foliation by complete isometric totally geodesic 
hypersurfaces called \emph{slices}. 

Fix a slice $M^n=M^n\times\{0\}$ and denote by $\Psi\colon\R\times M^n\to N^{n+1}$ 
the flux generated by $Y$.   The \emph{Killing graph} $\Sigma(u)$ 
associated to a function $u\in C^2(M)$ is the  hypersurface 
$$
\Sigma(u)=\{\Psi(u(x),x): x\in M^n\}.
$$
Recently we proved in \cite{DL} that under certain conditions any entire bounded 
Killing graph with constant mean curvature must be a slice.  Bounded means that 
the graph lies inside a slab, i.e., a set $M^n\times\mathbb{I}$ where 
$\mathbb{I}\subset\R$ is a closed interval.

It is a natural question if the condition that the graph is bounded can be 
weakened. In fact, this was recently shown to be the case by Rosenberg, Schulze 
and Spruck \cite{RSS}  for minimal graphs in a Riemannian product in $M^n\times\R$ 
where only a lower bound was required.  In turn, their result extended a well 
known theorem due to Bombieri, De Giorgi and Miranda \cite{BGM} for hypersurfaces 
in the Euclidean space that also holds for constant mean curvature due to 
independent work of Chern \cite{Ch} and Flanders \cite{Fl}.
\vspace{1,5ex}

 We denote by  $d(x)={\rm dist}(x,p)$ the geodesic distance from 
$p\in M^n$ and by $B_R(p)$  the geodesic ball centered at $p$ 
with radius $R>0$.  We also will make use of the 
notations $\rho_0=\sup_M\rho$ and  $\rho_1=\sup_M|\nabla\rho|$.

\begin{theorem}\po \label{main} 
Let $N^{n+1}=M^n\times_\rho\R$ be a Riemannian warped product manifold where 
$M^n$ is complete, noncompact such that $K_M\geq -K_0$ 
for $K_0\geq 0$ and $\mathrm{Ric}_M\geq 0$. Assume  $\inf_M\rho>0$ and 
$\|\rho\|_{C^2(M)}<\infty$. Then, any Killing graph with 
constant mean curvature $H\leq 0$ lying inside a region of the form
$$
\mathcal{R} = \left\{\Psi(s,x): 0\leq s\leq\frac{1}{\alpha\beta}
(C-\log \alpha\rho(x))\right\}
$$ 
for constants $\alpha>1$, $C>\log\alpha\rho_0$ and 
$\beta\geq n|H|\rho_0+2\rho_1$ positive must be a slice. 
\end{theorem}

From the above theorem one can easily recover the result in \cite{DL}
by a slight addition to the proof; see Remark \ref{remark} for details.
\vspace{1ex}

The proof of Theorem \ref{main} is done assuming 
$\|\rho\|_{C^1(M)}<\infty$ and ${\rm Ric}_N\geq-L$ for some constant $L\geq 0$.
That these conditions are weaker than  $\|\rho\|_{C^2(M)}<\infty$ 
follows from the relation between the   Ricci curvatures of $N^{n+1}$ and 
$M^n$ given by 
$$
\textrm{Ric}_N({\sf v}, {\sf v}) =\textrm{Ric}_M(\pi_*{\sf v}, \pi_*{\sf v})
-\frac{1}{\rho}\textrm{Hess}^M\rho(\pi_*{\sf v},\pi_*{\sf v})
-\frac{1}{\rho^3}\<{\sf v},Y\>^2\Delta^M\rho 
$$
where $\pi$ denotes the projection from $N^{n+1}$ to the factor $M^n$.
\vspace{1,5ex}

The region $\mathcal{R}$ in the above result is determined by the warping function
$\rho$ and may be unbounded as shown in the following example.

\begin{example}\po 
Let $M^n$ be a complete noncompact manifold with a pole at $p\in M^n$ and denote 
the radial coordinate  by $r$.  Take $\rho(r) = c+e^{-\psi(r)}$ for some constant 
$c>0$ and a smooth function $\psi>0$ with $\psi'>0$ such that $\psi(r)\to+\infty$ 
and $\psi'(r)\to 0$ as $r\to +\infty$.  In this situation,  the region 
$\mathcal{R}$ is unbounded.
\end{example}

The main achievement in this paper is a  strong improvement of the gradient 
estimate for Killing graphs when compared with the one in \cite{DL}. 
Here as there, we make use of the Korevaar-Simon method \cite{Ko} 
to show the existence of an a priori gradient estimate for a nonnegative solution  
of the corresponding PDE for constant mean curvature $H$ over a geodesic ball $B_R(p)$.  
\vspace{1,5ex}

 The function  $G\in C^\infty([0,+\infty))$ in the following result verifies 
that $f''=Gf$ where $f\in C^\infty([0,+\infty))$ satisfies $f(0)=0$, $f'(0)=1$, 
that $f>0$ outside the origin and $f'\geq 0$.  We also use the notation 
$\varrho_0=\sup_{B_R(p)}\rho$ and  $\varrho_1=\sup_{B_R(p)}|\nabla\rho|$.

\begin{theorem}\label{est}\po
Let $N^{n+1}=M^n\times_\rho\R$ be a complete Riemannian manifold 
satisfying  that ${\rm Ric}_N\geq-L$ for some constant $L\geq 0$.
Assume that the radial sectional curvatures on $M^n$ along the geodesics 
issuing from a fixed point $p\in M^n$  satisfy  $K^M_{\rm rad}\geq- G(d)$. 
Let $\Sigma(u)$ be a Killing graph with  constant mean curvature $H$ over 
$B_R(p)$ such that
\be\label{region}
0\leq u\leq\frac{1}{\alpha\beta}(C-\log\alpha\rho)
\ee
for constants $\alpha>1$, $C>\log\alpha\varrho_0$ and 
$\beta\geq n|H|\varrho_0+2\varrho_1$ positive.
Then, we have that
$$
|\nabla^M u(p)|\leq D
$$
where the constant $D=D(\alpha,C,\beta,|H|,u(p))$ is given by (\ref{c00}).
\end{theorem}

\proof It was shown in \cite{DHL} that a Killing graph $\Sigma(u)$ has mean 
curvature $H$ if and only if  the function $u\in C^2(M)$ satisfies the elliptic 
PDE of divergence form
$$
\textrm{div}_M\bigg(\frac{\nabla^M u}{W}\bigg)
-\frac{1}{2\gamma W}\<\nabla^M\gamma,\nabla^M u\>=nH\;\;
\mbox{where}\;\;\gamma=1/\rho^2.
$$
Here $H$ is computed with respect to the orientation of the Gauss map given by 
\be\label{gauss}
\mathcal{N}=\frac{1}{W}(\gamma Y-\Psi_*\nabla^M u)\;\;\mbox{where}
\;\;W=\sqrt{\gamma+|\nabla^M u|^2}.
\ee

In the sequel  $\nabla$ and $\Delta$ denote  the gradient and Laplace operator 
on $\Sigma(u)$.  If $H$ is constant it is  well-known \cite{FR} that
$$
\Delta\< Y,\mathcal{N}\> 
= -(|A|^2+\textrm{Ric}_N(\mathcal{N}))\<Y,\mathcal{N}\>
$$
where $|A|$ stands for the norm of the second fundamental form of  $\Sigma(u)$. 
Thus,  using that  $\<Y,\mathcal{N}\> = 1/W$ we obtain
\be\label{four}
\Delta W-\frac{2}{W}|\nabla W|^2=(|A|^2+\textrm{Ric}_N(\mathcal{N}))W.
\ee

We consider the functions 
$$
h(d)=\frac{1}{C_R}\int_0^d f(\tau)d\tau\;\;\;\mbox{where}\;\;\;
C_R=\int_0^R f(\tau)d\tau
$$
and
$$
\xi(s)=e^C\int_0^s e^{-\alpha\beta\tau}\,{\rm d}\tau
$$
where $s$ parametrizes the line factor $\R$ in $N^{n+1}$. 
Let $\mathcal{K}^+$ denote the solid half-cylinder $\Psi(\R^+\times B_R)$
with $B_R=B_R(p)$ and let  $\phi\colon\mathcal{K}^+\to\R$ be the function
$$
\phi(y)=\left(1-h(d(x))-C_0\xi(s)\right)^+
\;\;\textrm{if}\;\;  y=\Psi(s,x)\in\mathcal{K}^+
$$
where $C_0=1/2\,\xi(u(p))>0$ and $+$ means taking the positive part. 
Let $\eta$ be the function supported in the portion of the graph $\Sigma (u)$ 
over $B_R$ given by
$$
\eta(y) = e^{K\phi(y)}-1, \quad y=\Psi(u(x),x),
$$
for some constant $K>0$ and set $U=\eta W$. 

Let $\Psi(u(q),q)$ with $q\in B_R$ be 
a (necessarily interior) point of maximum of $U$.
We first consider the case when $q\in B_R\backslash C(p)$, where $C(p)$ 
is the cut locus of $p$ in $M^n$, thus $U$ is smooth near $q$.  
Without further reference in the sequel we compute at $\Psi(u(q),q)$.  
It holds that
$$
\nabla U=\eta\nabla W+W\nabla\eta=0
$$
and
$$
\Delta U=W\Delta\eta+\Big(\Delta W-\frac{2}{W}|\nabla W|^2\Big)\eta.
$$
Since the Hessian form of $U$ is nonpositive, 
we have from (\ref{four}) that
$$
\Delta U=W\big(\Delta\eta +(|A|^2
+\textrm{Ric}_N(\mathcal{N}))\eta\big)\leq 0.
$$
Making use of the Ricci curvature assumption, we obtain 
$$
\Delta\eta\leq L\eta
$$
or, equivalently, that 
\be\label{zz}
\Delta\phi + K|\nabla \phi|^2\leq\frac{L}{K}\cdot
\ee

In the sequel, we estimate both terms in the left hand side of (\ref{zz}). 
We have $\bar\nabla s=\gamma Y$  where $\bar\nabla$ denotes the gradient in $N^{n+1}$
but also will stand for the Riemannian connection  in $N^{n+1}$.
Where the function $\phi$ is differentiable and positive, we have
\begin{eqnarray}
 \bar\nabla \phi (y)\!\!\!&=&\!\!\!
-h'(d)\bar\nabla d (y)-C_0\dot\xi\gamma Y(y)\nonumber\\
\!\!\!&=&\!\!\!-h'(d)\Psi_*\nabla^Md(x)-C_0\dot\xi\gamma Y(y)\label{gradphi}
\end{eqnarray}
where $h'=h'(d)$ and $\cdot$ indicates the derivative with respect to $s$. Then,
\be\label{old}
\<\bar\nabla\phi,\mathcal{N}\>=\frac{h'}{W}\<\nabla^Md,\nabla^Mu\>
-\frac{\gamma}{W}C_0\dot\xi.
\ee
Therefore, 
\bea
|\nabla\phi|^2\!\!\!&=&\!\!\! |\bar\nabla \phi|^2
-\<\bar\nabla\phi,\mathcal{N}\>^2\\
\!\!\!&=&\!\!\! h'^2\Big(1-\frac{1}{W^2}\<\nabla^Md,\nabla^Mu\>^2\Big) 
+\frac{\gamma C_0^2\dot\xi^2}{W^2}|\nabla^Mu|^2
+\frac{2\gamma h'C_0\dot\xi}{W^2}\<\nabla^Md,\nabla^Mu\>\\
\!\!\!&\geq &\!\!\!h'^2\Big(1-\frac{1}{W^2}|\nabla^M u|^2\Big)
+\frac{\gamma C_0^2\dot\xi^2}{W^2}|\nabla^Mu|^2
-\frac{2\gamma h'C_0\dot\xi}{W^2}|\nabla^M u|\\
\!\!\!&=&\!\!\!\gamma\dot\xi^2\Big(\frac{C_0}{W}|\nabla^M u|
-\frac{h'}{\dot\xi W}\Big)^2.
\eea
It follows from the assumption (\ref{region})  that 
$\dot\xi(s)\geq\alpha\rho(x)$, i.e.,
$\sqrt\gamma\,\dot\xi\geq\alpha$.
Hence,
$$
\left(\frac{C_0}{W}|\nabla^M u|
-\frac{h'}{\dot\xi W}\right)^2\geq\frac{C_0^2}{\gamma\dot\xi^2}
$$
is implied by 
$$
\left(\frac{C_0}{W}|\nabla^M u|
-\frac{h'}{\dot\xi W}\right)^2\geq\frac{C_0^2}{\alpha^2}
$$
that, in turn, is equivalent to
$$
\left(\frac{C_0}{W}|\nabla^M u|-\frac{h'}{\dot\xi W}
-\frac{C_0}{\alpha}\right)\left(\frac{C_0}{W}|\nabla^M u|
-\frac{h'}{\dot\xi W}+\frac{C_0}{\alpha}\right)\geq 0.
$$
Clearly, the latter holds if the first factor is nonnegative
or, equivalently, if
\be\label{ineq-main}
|\nabla^M u|-\frac{1}{\alpha} W\ge\frac{h'}{\dot\xi C_0}\cdot
\ee
Assume that 
\be\label{delta}
\frac{|\nabla^M u|^2}{\gamma}\geq\frac{1}{\alpha^2-1}\cdot
\ee
It follows that 
$$
|\nabla^M u|-\frac{1}{\alpha} W\geq 0
$$
and thus (\ref{ineq-main}) is implied by
\be\label{ineq-main-2}
|\nabla^M u|^2-\frac{2}{\alpha}|\nabla^M u|W 
+\frac{1}{\alpha^2}W^2\geq\frac{h'^2}{\dot\xi^2C_0^2}\cdot
\ee
Using that
$$
2|\nabla^M u|W\leq |\nabla^M u|^2+W^2
$$
it follows that (\ref{ineq-main-2}) is implied by 
$$
|\nabla^M u|^2-\frac{1}{\alpha}(|\nabla^M u|^2+W^2) 
+\frac{1}{\alpha^2}W^2\geq\frac{h'^2}{\dot\xi^2C_0^2}\cdot
$$
But this is equivalent to 
$$
\left(1-\frac{1}{\alpha}\right)^2 \frac{|\nabla^M u|^2}{\gamma}
-\frac{1}{\alpha}\left(1-\frac{1}{\alpha}\right)\geq\frac{h'^2}{\gamma\dot\xi^2 C_0^2}
$$
that is implied by
$$
\left(1-\frac{1}{\alpha}\right)^2\frac{|\nabla^M u|^2}{\gamma}
-\frac{1}{\alpha}\left(1-\frac{1}{\alpha}\right)\geq\frac{h'^2}{\alpha^2C_0^2}
$$
Now assume that 
\be\label{hip-grad-2}
\frac{|\nabla^Mu|^2}{\gamma}
\geq\frac{1}{\alpha-1}+\frac{h'^2}{(\alpha-1)^2C_0^2}\cdot
\ee  
Since (\ref{hip-grad-2}) implies (\ref{delta}), we conclude that the estimate
\be\label{first}
|\nabla\phi|^2\geq C_0^2.
\ee
for the first term in (\ref{zz}) holds under (\ref{hip-grad-2}).
\vspace{1ex}

We now estimate the second term in (\ref{zz}).
Taking a local orthonormal tangent frame $\{e_i\}_{i=1}^n$ in $\Sigma (u)$, we have 
\begin{eqnarray}\label{one}
\Delta \phi\!\!\!&=&\!\!\!\sum_{i=1}^n\< \nabla_{e_i}\nabla\phi, e_i\>\nonumber\\ 
\!\!\!&=&\!\!\!\sum_{i=1}^n\< \bar\nabla_{e_i}(\bar\nabla\phi
-\<\bar\nabla\phi,\mathcal{N}\>\mathcal{N}), e_i\>\nonumber\\
\!\!\!&=&\!\!\!\sum_{i=1}^n \<\bar\nabla_{e_i}\bar\nabla \phi,e_i\>
+\<\bar\nabla\phi,\mathcal{N}\>nH.
\end{eqnarray}
We obtain from (\ref{gradphi}) that
\be\label{new}
\<\bar\nabla_{e_i}\bar\nabla\phi,e_i\>=-h'\<\bar\nabla_{e_i}
\bar\nabla d,e_i\>-h''\<\bar\nabla d,e_i\>^2
-C_0\<\bar\nabla_{e_i}\dot\xi\bar\nabla s,e_i\>.
\ee
On the other hand,
$$
\<\bar\nabla_{e_i}\bar\nabla d,e_i\>=\< \nabla^M_{\hat{e}_i}\nabla^M d,\hat{e}_i\>
+\gamma^2\<Y,e_i\>^2\<\bar\nabla_Y\bar\nabla d,Y\>
$$
where $e_i=\hat{e}_i+\<Y,e_i\>\gamma Y$.
 From our assumptions and the comparison theorem for the Hessian 
(see Theorem 2.3 in \cite{PRS}) we obtain
$$
\nabla^M \nabla^M d \leq g(d)(\<\cdot,\cdot\>-\textrm{d}d\otimes\textrm{d}d)
$$
where $g=f'/f$.  And since the metric is a warped product, we have
$$
\bar\nabla_Y\bar\nabla d= \frac{1}{\rho}\<\nabla^Md,\nabla^M\rho\>Y.
$$
It follows  that
$$
\<\bar\nabla_{e_i}\bar\nabla d,e_i\>\leq 
g(\<\hat e_i,\hat e_i\>-\<\nabla^M d,\hat e_i\>^2)
+\frac{1}{2}\<Y,e_i\>^2|\nabla^M\gamma|.
$$
Since $h'\geq 0$, we have from (\ref{new}) that
\bea
\<\bar\nabla_{e_i}\bar\nabla \phi, e_i\> 
\!\!\!&\geq&\!\!\! -h'(g(1-\gamma\<Y,e_i\>^2-\<\nabla^M d,\hat{e}_i\>^2)
+\frac{1}{2}\<Y,e_i\>^2|\nabla^M\gamma|)\\
\!\!\!&&\!\!\!-h''\<\nabla^M d,\hat{e}_i\>^2
-C_0\<\bar\nabla_{e_i}\dot\xi\bar\nabla s,e_i\>.
\eea
On one hand,
\bea
\sum_{i=1}^n\<\bar\nabla_{e_i}\dot\xi\bar\nabla s,e_i\>
\!\!\!&=&\!\!\!\sum_{i=1}^n\dot\xi\<\bar\nabla_{e_i}\bar\nabla s,e_i\>
+\ddot\xi\sum_{i=1}^n\<\bar\nabla s,e_i\>^2\\
\!\!\!&=&\!\!\!\dot\xi(\Delta s-\<\bar\nabla s,\mathcal{N}\>nH)
+\ddot\xi|\nabla s|^2.
\eea
On the other hand,
\bea
\Delta s\!\!\!&=&\!\!\!\sum_{i=1}^n\<\nabla_{e_i}\gamma Y^T,e_i\> \\
\!\!\!&=&\!\!\!\< \nabla \gamma,Y^T\> 
+\gamma\sum_{i=1}^n\<\bar\nabla_{e_i}Y,e_i\>+n H\<\gamma Y,\mathcal{N}\>\\
\!\!\!&=&\!\!\!\<\bar\nabla\gamma,Y^T\>
+\<\bar\nabla s,\mathcal{N}\>nH
\eea
where $Y^T$ denotes the component of $Y$ tangent to $\Sigma({u})$.
Hence,
$$
\sum_{i=1}^n\<\bar\nabla_{e_i}\dot\xi\bar\nabla s,e_i\> 
=\dot\xi\<\bar\nabla\gamma,Y^T\>+\ddot\xi|\nabla s|^2.
$$
Using that $h''=h'g$, we obtain 
$$
\sum_{i=1}^n\<\bar\nabla_{e_i}\bar\nabla\phi,e_i\>\geq -
h'(g(n-\gamma|Y^T|^2)+(1/2)|\nabla^M\gamma||Y^T|^2)
-C_0\dot\xi\<\bar\nabla\gamma,Y^T\>
-C_0\ddot\xi|\gamma Y^T|^2. 
$$
We have from (\ref{gauss}) and (\ref{old}) that
$$
\<\bar\nabla\phi,\mathcal{N}\>=
-h'\<\bar\nabla d,\mathcal{N}\>-\frac{\gamma}{W}C_0\dot\xi.
$$
Setting 
$$
\kappa=\frac{1}{\rho}|\nabla^M\rho|=\frac{1}{2\gamma}|\nabla^M\gamma|
$$
we obtain from (\ref{one}) that
$$
\Delta\phi\geq-h'(ng+\kappa\gamma|Y^T|^2+n|H|)
-\frac{\gamma}{W}C_0\dot\xi nH-C_0\dot\xi\<\bar\nabla\gamma,Y^T\>
-C_0\ddot\xi|\gamma Y^T|^2.
$$
Hence,
$$
\Delta\phi\geq-\frac{n}{C_R}f'(d)
-(\kappa\gamma |Y^T|^2+n|H|)\frac{1}{C_R}f(d)
-C_0\left(\frac{\gamma}{W}\dot\xi nH
+\dot\xi\<\bar\nabla\gamma,Y^T\>+\ddot\xi |\gamma Y^T|^2\right).
$$
Since $Y^T=Y-(1/W)\mathcal{N}$, we have
\be\label{rel}
|\gamma Y^T|^2=\frac{\gamma}{W^2}|\nabla^M u|^2.
\ee
Thus,
$$
\frac{\gamma}{W}=|\gamma Y^T|\frac{\sqrt\gamma}{|\nabla^M u|}.
$$
Assume that
\be\label{hyp-grad-3}
|\nabla^M u|^2\geq\gamma.
\ee
Since $\dot\xi\ge 0$ and
$$
\<\bar\nabla\gamma,Y^T\>\leq 2\kappa|\gamma Y^T|,
$$
then we have using (\ref{hyp-grad-3}) that
$$
\frac{\gamma}{W}\dot\xi nH+\dot\xi\<\bar\nabla\gamma,Y^T\>
+\ddot\xi|\gamma Y^T|^2\leq|\gamma Y^T|(\dot\xi(n|H|
+2\kappa)+\ddot\xi|\gamma Y^T|).
$$
 From (\ref{delta}) and (\ref{rel}) we have
$$
|\gamma Y^T|\geq\frac{1}{\alpha}\sqrt\gamma.
$$
Considering that $\ddot\xi=-\alpha\beta\dot\xi\leq 0$ we obtain that
$$
\frac{\gamma}{W}\dot\xi nH+\dot\xi\<\bar\nabla\gamma,Y^T\>
+\ddot\xi|\gamma Y^T|^2
\leq|\gamma Y^T|\sqrt\gamma\dot\xi(n|H|\rho+2|\nabla^M\rho|-\beta).
$$
Now the choice of $\beta$ gives that 
$$ 
\Delta\phi\geq-\frac{n}{C_R}f'(d)
-(\kappa\gamma|Y^T|^2+n|H|)\frac{1}{C_R}f(d).
$$
We have from (\ref{rel}) that $\gamma|Y^T|\leq 1$ and from $f'\geq 0$ 
that $f(d)\leq f(R)$.  Hence,
$$
\Delta\phi\geq-n\frac{f'(d)}{C_R}
-(\kappa+n|H|)\frac{f(R)}{C_R}.
$$
It follows that
\be\label{est2}
\Delta\phi>-2A
\ee
where 
$$
2A=\frac{n}{C_R}\sup_{B_R} f'(d)
+\frac{f(R)}{C_R}(n|H|+\sup_{B_R}\kappa). 
$$
 From (\ref{zz}), (\ref{first}) and (\ref{est2}) we obtain 
\be\label{zzz}
\frac{L}{K}\geq \Delta\phi + K|\nabla \phi|^2>-2A+KC_0^2.
\ee
Taking
\be\label{k0}
K>\frac{1}{C_0^2}\left(A+\sqrt{A^2+C_0^2L}\right),
\ee
we obtain a contradiction in (\ref{zzz}). 
Thus,  we conclude from (\ref{hip-grad-2}) and (\ref{hyp-grad-3})  that 
\be\label{c0}
\frac{|\nabla^M u|^2}{\gamma}(q)< D_0
=\max_{B_R}\left\{1,\frac{1}{\alpha-1}
+\frac{4f^2(R)\xi^2(u(p))}{(\alpha-1)^2C_R^2}\right\}.
\ee
Then,
$$
W(q)\leq D_1=\frac{1}{r_0}\sqrt{1+D_0}
$$
where $r_0=\min_{B_R}\rho$. Hence,
$$
U(p)=(e^{K/2}-1)W(p)\leq U(q)\leq D_1(e^K-1)
$$
where $K$ is given by (\ref{k0}).
It follows that
\be\label{c00}
|\nabla^M u(p)|\leq D=D_1(e^{K/2}+1).
\ee

To conclude the proof we observe that the same argument given in \cite{RSS}
proves that the restriction that $q\not\in C(p)$ can be dropped. \hfill $\square$
\vspace{2ex}

\noindent \emph{Proof of Theorem \ref{main}:}  
We claim that there is a global gradient estimate obtained by taking $R\to +\infty$ 
in Theorem \ref{est}. To see this we choose the function 
$$
f(t)=\frac{1}{\sqrt{K_0}}\sinh\sqrt{K_0}\,t
$$
that satisfies all the requirements in that result. 
Then, we have $G=K_0$ and 
$$
C_R=\frac{1}{K_0}\left(\cosh\sqrt{K_0}R-1\right).
$$
Since $f(R)/C(R)\to 1$ as $R\to +\infty$, it follows from (\ref{c0}) that we can take 
$$
D_0=\max\left\{1,\frac{1}{\alpha-1}+\frac{4K_0\xi^2(u(p))}{(\alpha-1)^2}\right\}.
$$
We also have that
$$
2A\leq 2A_1=nK_0+(n|H|+\rho_1/r_0)\sqrt{K_0}
$$
where $r_0=\inf_M\rho$. Since $0\leq\xi<e^C/\alpha\beta$, then
$$
C_0\geq C_1=\alpha\beta/2e^C.
$$
Thus, from (\ref{k0}) we can take
$$
K>\frac{1}{C_1^2}\left(A_1+\sqrt{A_1^2+C_1^2L}\right),
$$
and the claim follows from (\ref{c00}).  \newpage

We now make use of part of the argument in \cite{DL} to prove that $H=0$.  
It was shown there the Omori-Yau maximum principle for the Laplacian
holds on $\Sigma(u)$.  Moreover, we have that $u=s|_\Sigma$ 
satisfies
\be\label{final}
\Delta u =\frac{n\gamma}{\sqrt{\gamma+|\nabla^M u|^2}}H
+\frac{1}{\gamma}\<\nabla\gamma,\nabla u\>.
\ee
Being $u$ a bounded function from below on $\Sigma(u)$, the  Omori-Yau
maximum principle assures that there exists a sequence 
$\{y_k\}_{k\in\mathbb{N}}$ such that
$$
|\nabla u(y_k)|<1/k\;\;\; \mbox{and}\;\;\;\;\Delta u(y_k)>-1/k.
$$
Since we have already shown that the coefficient of $H$ is bounded from 
below by a positive number,  we obtain from (\ref{final}) that $H=0$.

\begin{remark}\po\label{remark} {\em Notice that the above argument
can be used to show that a graph of any constant mean curvature $H$ inside 
a slab must be minimal. 
}\end{remark}

We conclude the proof with an argument from \cite{DL} or \cite{RSS}.
Now the PDE can be written as $Lu=0$ where
$$
Lu=e^{-\varphi}\textrm{div}_\Sigma e^{\varphi}\nabla u,\;\;\;\varphi=2\log\rho.
$$
In a system of coordinates $\{x^i\}_{i=1}^n$ 
in $M^n$ with $\sigma_{ij}=\<\partial_{x^i},\partial_{x^j}\>$, we have
$$
L=e^{-\varphi}\textrm{div}_\Sigma (e^\varphi g^{ij}u_i\partial_{x^j})
$$
where 
$$
g^{ij}=\sigma^{ij}-\frac{u^iu^j}{W^2}\;\;\mbox{and}\;\;u^i=\sigma^{ik}u_k.
$$
Since there is a global gradient estimate, then $L$ is an uniformly 
elliptic operator in divergence form. Hence, if we view $L$ as a operator 
acting on $M^n$ and since $\textrm{Ric}_M\geq 0$, it follows from 
Theorem 7.4 in \cite{SC} that $u$ is constant. \hfill $\square$

{\renewcommand{\baselinestretch}{1} \hspace*{-20ex}\begin{tabbing}
\indent \= Marcos Dajczer\\
\> IMPA \\
\> Estrada Dona Castorina, 110\\
\> 22460-320 -- Rio de Janeiro -- Brazil\\
\> marcos@impa.br\\
\end{tabbing}}

\vspace*{-4ex}

{\renewcommand{\baselinestretch}{1} \hspace*{-20ex}\begin{tabbing}
\indent \= Jorge Herbert S. de Lira\\
\> Departamento de Matematica - UFC, \\
\> Bloco 914 -- Campus do Pici\\
\> 60455-760 -- Fortaleza -- Ceara -- Brazil\\
\> jorge.lira@mat.ufc.br
\end{tabbing}}

\end{document}